\magnification 1200

%
%
\newdimen\FigSize	\FigSize=.9\hsize 
%
\newskip\abovefigskip	\newskip\belowfigskip
\gdef\epsfig#1;#2;{\par\vskip\abovefigskip\penalty -500
   {\everypar={}\epsfxsize=#1\noindent
    \centerline{\epsfbox{#2}}}%
    \vskip\belowfigskip}%
%
\newskip\figtitleskip
\gdef\tepsfig#1;#2;#3{\par\vskip\abovefigskip\penalty -500
   {\everypar={}\epsfxsize=#1\noindent
    \vbox
      {\centerline{\epsfbox{#2}}\vskip\figtitleskip
       \centerline{\figtitlefont#3}}}%
    \vskip\belowfigskip}%
%
\newcount\FigNr	\global\FigNr=0
\gdef\nepsfig#1;#2;#3{\global\advance\FigNr by 1
   \tepsfig#1;#2;{Figure\space\the\FigNr.\space#3}}%
%
%
%
\gdef\ipsfig#1;#2;{
   \midinsert{\everypar={}\epsfxsize=#1\noindent
	      \centerline{\epsfbox{#2}}}%
   \endinsert}%
%
\gdef\tipsfig#1;#2;#3{\midinsert
   {\everypar={}\epsfxsize=#1\noindent
    \vbox{\centerline{\epsfbox{#2}}%
          \vskip\figtitleskip
          \centerline{\figtitlefont#3}}}\endinsert}%
%
\gdef\nipsfig#1;#2;#3{\global\advance\FigNr by1%
  \tipsfig#1;#2;{Figure\space\the\FigNr.\space#3}}%
\newread\epsffilein    
\newif\ifepsffileok    
\newif\ifepsfbbfound   
\newif\ifepsfverbose   
\newdimen\epsfxsize    
\newdimen\epsfysize    
\newdimen\epsftsize    
\newdimen\epsfrsize    
\newdimen\epsftmp      
\newdimen\pspoints     
\pspoints=1bp          
\epsfxsize=0pt         
\epsfysize=0pt         
\def\epsfbox#1{\global\def\epsfllx{72}\global\def\epsflly{72}%
   \global\def\epsfurx{540}\global\def\epsfury{720}%
   \def\lbracket{[}\def\testit{#1}\ifx\testit\lbracket
   \let\next=\epsfgetlitbb\else\let\next=\epsfnormal\fi\next{#1}}%
\def\epsfgetlitbb#1#2 #3 #4 #5]#6{\epsfgrab #2 #3 #4 #5 .\\%
   \epsfsetgraph{#6}}%
\def\epsfnormal#1{\epsfgetbb{#1}\epsfsetgraph{#1}}%
\def\epsfgetbb#1{%
%
%
\openin\epsffilein=#1
\ifeof\epsffilein\errmessage{I couldn't open #1, will ignore it}\else
%
%
   {\epsffileoktrue \chardef\other=12
    \def\do##1{\catcode`##1=\other}\dospecials \catcode`\ =10
    \loop
       \read\epsffilein to \epsffileline
       \ifeof\epsffilein\epsffileokfalse\else
%
%
          \expandafter\epsfaux\epsffileline:. \\%
       \fi
   \ifepsffileok\repeat
   \ifepsfbbfound\else
    \ifepsfverbose\message{No bounding box comment in #1; using defaults}\fi\fi
   }\closein\epsffilein\fi}%
%
%
\def\epsfsetgraph#1{%
   \epsfrsize=\epsfury\pspoints
   \advance\epsfrsize by-\epsflly\pspoints
   \epsftsize=\epsfurx\pspoints
   \advance\epsftsize by-\epsfllx\pspoints
%
%
   \epsfxsize\epsfsize\epsftsize\epsfrsize
   \ifnum\epsfxsize=0 \ifnum\epsfysize=0
      \epsfxsize=\epsftsize \epsfysize=\epsfrsize
%
%
     \else\epsftmp=\epsftsize \divide\epsftmp\epsfrsize
       \epsfxsize=\epsfysize \multiply\epsfxsize\epsftmp
       \multiply\epsftmp\epsfrsize \advance\epsftsize-\epsftmp
       \epsftmp=\epsfysize
       \loop \advance\epsftsize\epsftsize \divide\epsftmp 2
       \ifnum\epsftmp>0
          \ifnum\epsftsize<\epsfrsize\else
             \advance\epsftsize-\epsfrsize \advance\epsfxsize\epsftmp \fi
       \repeat
     \fi
   \else\epsftmp=\epsfrsize \divide\epsftmp\epsftsize
     \epsfysize=\epsfxsize \multiply\epsfysize\epsftmp   
     \multiply\epsftmp\epsftsize \advance\epsfrsize-\epsftmp
     \epsftmp=\epsfxsize
     \loop \advance\epsfrsize\epsfrsize \divide\epsftmp 2
     \ifnum\epsftmp>0
        \ifnum\epsfrsize<\epsftsize\else
           \advance\epsfrsize-\epsftsize \advance\epsfysize\epsftmp \fi
     \repeat     
   \fi
%
%
   \ifepsfverbose\message{#1: width=\the\epsfxsize, height=\the\epsfysize}\fi
   \epsftmp=10\epsfxsize \divide\epsftmp\pspoints
   \vbox to\epsfysize{\vfil\hbox to\epsfxsize{%
      \includegraphics{#1}%
      \hfil}}%
\epsfxsize=0pt\epsfysize=0pt}%
%
%
{\catcode`\%=12 \global\let\epsfpercent=
%
%
\long\def\epsfaux#1#2:#3\\{\ifx#1\epsfpercent
   \def\testit{#2}\ifx\testit\epsfbblit
      \epsfgrab #3 . . . \\%
      \epsffileokfalse
      \global\epsfbbfoundtrue
   \fi\else\ifx#1\par\else\epsffileokfalse\fi\fi}%
%
%
\def\epsfgrab #1 #2 #3 #4 #5\\{%
   \global\def\epsfllx{#1}\ifx\epsfllx\empty
      \epsfgrab #2 #3 #4 #5 .\\\else
   \global\def\epsflly{#2}%
   \global\def\epsfurx{#3}\global\def\epsfury{#4}\fi}%
%
%
\def\epsfsize#1#2{\epsfxsize}%
%
%

\epsfverbosetrue			
\abovefigskip=\baselineskip		
\belowfigskip=\baselineskip		
\global\let\figtitlefont\bf		
\global\figtitleskip=.5\baselineskip	

\let\nd\noindent %
\font\tenmsb=msbm10
\font\sevenmsb=msbm7
\font\fivemsb=msbm5
\newfam\msbfam
\textfont\msbfam=\tenmsb
\scriptfont\msbfam=\sevenmsb
\scriptscriptfont\msbfam=\fivemsb
\def\Bbb#1{\fam\msbfam\relax#1}

\def\Z{{\Bbb Z}}

\def\I{{\Bbb I}}

\def\a{\alpha}
\def\b{\beta}
\def\qed{\hbox{\hskip 6pt\vrule width6pt height7pt depth1pt \hskip1pt}}
\def\vs#1 {\vskip#1truein}
\def\hs#1 {\hskip#1truein}

\def\Month{\ifcase\number\month \relax\or January \or February \or
  March \or April \or May \or June \or July \or August \or September
  \or October \or November \or December \else \relax\fi }
\def\date{\Month \the\day, \the\year}

  \hsize=6truein	\hoffset=.25truein 
  \vsize=8.8truein	
  \pageno=1     \baselineskip=12pt
  \parskip=3 pt		\parindent=20pt
  \overfullrule=0pt	\lineskip=0pt	\lineskiplimit=0pt
  \hbadness=10000 \vbadness=10000 
\pageno=0

\footline{\ifnum\pageno=0\hss\else\hss\tenrm\folio\hss\fi}
\hbox{}
\vskip 1truein\centerline{{\bf SUBGROUPS OF SO(3) ASSOCIATED WITH TILINGS}}
\vskip .5truein\centerline{by}
\centerline{Charles Radin${}^1$ and Lorenzo Sadun${}^2$}
\footnote{}{1\ Research supported in part by NSF Grant No. DMS-9304269 and
\vs-.1 \hs.15 Texas ARP Grant 003658-113\hfil}
\footnote{}{2\ Research supported in part by an NSF Mathematical Sciences 
Postdoctoral 
\vs-.1 \hs.15 Fellowship and Texas ARP Grant 003658-037 \hfil}
\vskip .2truein\centerline{Mathematics Department}
\centerline{University of Texas}
\centerline{Austin, TX\ \ 78712}
\vs.1
\centerline{radin@math.utexas.edu and sadun@math.utexas.edu}
\vs.5
\centerline{{\bf Abstract}}
\vs.1 \nd
We give a thorough analysis of those subgroups of $SO(3)$ generated by
rotations about perpendicular axes by $2\pi/p$ and $2\pi/q$. A
corollary is that such a group is the free product of the cyclic
groups of rotations about the separate axes if and only if $p,q\ge 3$
and are both odd. These groups are naturally associated with a family
of hierarchical tilings of Euclidean 3-space.
\vs.8
\vs.2
\centerline{Subject Classification:\ \ 51F25, 52C22}
\vfill\eject 

\nd {\bf \S 0.\ Introduction}
\vs.1
We analyze certain subgroups of $SO(3)$ motivated by polyhedral
tilings of Euclidean 3-space. The tilings are made by a general
iterative procedure and consist of congruent copies of some finite set
of polyhedra. Our interest here is in the relative orientations of the
polyhedra in such a tiling. 

An example called ``quaquaversal'' tilings, consisting of congruent
copies of a single triangular prism, was analyzed in [CoR]. Let
$G(p,q)$ be the group of rotations in 3-space generated by rotations, 
about perpendicular axes, by $2\pi/p$ and $2\pi/q$. In a quaquaversal 
tiling the orientations of any two prisms are related by an element of
$G(6,4)$. More precisely, in any cube of side $2^n$, the relative 
orientations of the prisms are words of length $j(n)$ or less
in the generators of $G(6,4)$, where $j(n)$ is a function that grows
linearly with $n$.  All words of length $k(n)$ or less arise in this
manner, where $k(n)$ also grows linearly with $n$. To
show that the number of distinct orientations in such a cube grows
exponentially in $n$, presentations of the subgroups $G(3,4)$ and
$G(3,3)$ were derived in [CoR]. In that paper the rotations are
represented by conjugation of real quaternions, and the results are
unique factorizations for certain elements in {\it noncommutative}
subrings of quaternions.

In this paper we extend the analysis to all groups $G(p,q)$, $p,q\ge
3$, using a somewhat simpler method. We represent the rotations by
explicit $SO(3)$ matrices, and obtain our factorization theorems by
applying {\it commutative} ring theory to the individual matrix
elements.  We prove that $G(p,q)$ is the free product of the cyclic
groups of rotations about the separate axes if and only if both $p$
and $q$ are odd.  If $p$ or $q$ is even, 
we can write $G(p,q)$ as the amalgamated free product
of two finite groups.  In all cases we obtain
canonical forms for the elements of $G(p,q)$. (We note the following:
if $p$ or $q$ is 1, $G(p,q)$ is a finite cyclic group; if $p$ or $q$
is 2, $G(p,q)$ is a finite dihedral group; $G(4,4)$ is the finite
group of symmetries of the cube; all the other $G(p,q)$ are dense in
$SO(3)$.)

The organization of the paper is as follows.  Our main results,
presentations of the groups $G(p,q)$ and canonical forms for the group
elements, are given in $\S 1$ and $\S 2$, respectively. In $\S 3$ we
describe a new tiling for which the orientations of the polyhedra are
given by $G(10,4)$. In $\S 4$ we consider orientation groups
$G(\theta,4)$, where the angle of rotation $\theta$ is an irrational
multiple of $2 \pi$.  First we describe a new tiling with orientation
group containing $G(\tan^{-1}({4 \over 3}),4)$ and find a presentation
of that group.  Then we discuss the case of $G(\omega,4)$ where
$\exp(i \omega)$ is transcendental.
\vs.3
\nd {\bf \S 1.\ Presentations for G(p,q)}
\vs.1
In this section we state and prove a classification theorem for the
groups $G(p,q)$.  After stating the theorem and deriving some corollaries,
we begin the proof.  We first reduce the theorem to Lemma 1, then 
reduce Lemma 1 to Lemmas 2 and 3, and then prove Lemmas 2 and 3.
Although the statement of the theorem is group theory, much of the proof,
and in particular Lemmas 2 and 3, is commutative ring theory.

Given positive integers $p,\ell$ and $q$, we define rotations
$A=R_x^{2\pi/p}$, $L=R_y^{2\pi/\ell}$, $S= R_y^{2\pi/4}$ 
and $B=R_z^{2\pi/q}$, where $R_x^{2\pi/p}$ is a rotation about the
$x$-axis by angle $2 \pi/p$, etc.  Let $G(p,\ell,q)$ be the group 
generated by $A,L$ and $B$, and let $G(p,q)\equiv
G(p,1,q)$ be the group generated by $A$ and $B$. 
\vs.1
\nd {\bf Theorem 1: Presentations for G(p,q)}

{\it \nd (i) If $p,q\ge 3$ are odd, then $G(p,q)$ is isomorphic to
the free product
$$ \Z_p * \Z_q =\ <\a,\b\,:\,\a^p,\,\b^q>.\eqno(1.1)$$

\nd (ii) If $p\ge 4$ is even and $q\ge 3$ is odd, then $G(p,q)$ has the
presentation 
$$<\a,\b\,:\,\a^p,\,\b^q,\,\a^{p/2}\b\a^{p/2}\b>. \eqno(1.2)$$

\nd (iii) If $p\ge 4$ is even and $q=2s$, $s\ge 3$ odd, then $G(p,q)$ has the
presentation 
$$<\a,\b\,:\,\a^p,\,\b^q,\,\a^{p/2}\b\a^{p/2}\b,\,\b^{q/2}\a\b^{q/2}\a>.
\eqno (1.3)$$

\nd (iv) If $4$ divides both $p$ and $q$, then $G(p,1,q)=G(lcm(p,q),4,1)$.

In cases (i), (ii) and (iii), the isomorphism between the abstract
presentation and $G(p,q)$ is given by $\a \mapsto A$, $\b \mapsto B$.}
\vs.1
These results can be rephrased in terms of free products and 
amalgamated free products.  

\nd {\bf Corollary 1:}\ {\it If $p,q\ge 3$, then $G(p,q)$ is 
isomorphic to the free product
$$ \Z_p * \Z_q =\ <\a,\b\,:\,\a^p,\,\b^q>, \eqno(1.4)
$$
\nd with the isomorphism given by $\a \mapsto A$, $\b \mapsto B$,
if and only if both $p$ and $q$ are odd.}
\vs.1
Now let $D_p$ denote the dihedral group $<\a, \gamma : \a^p, \gamma^2,
\gamma \a \gamma \a>$.  In cases (ii) and (iii) we can introduce a
new generator $\mu$, which we then set equal to $\a^{p/2}$, and in
case (iii) we introduce $\gamma$, which we set equal to $\b^{q/2}$.  The
subgroup of $G(p,q)$ generated by $\a$ and $\gamma$ is then $D_p$, while
the subgroup generated by $\b$ and $\mu$ is $D_q$.  In case 
(iii), $\gamma$ and $\mu$ generate a $D_2$ subgroup.

\nd {\bf Corollary 2:}\ {\it If $p\ge 4$ is even and $q\ge 3$ is odd, then 
$G(p,q)$ is isomorphic to the amalgamated free product
$$ \Z_p *\hs-.03 {\lower.5ex\hbox{}}_{\Z_2} D_q, \eqno (1.5) $$
where $\a^{p/2} \in \Z_p$ is identified with $\mu \in D_q$.

If $p\ge 4$ is even and $q=2s$, $s\ge 3$ odd, then $G(p,q)$
is isomorphic to the amalgamated free product
$$ D_p *\hs-.045 {\lower.5ex\hbox{}}_{D_2} D_q, \eqno (1.6) $$
where $\gamma \in D_p$ is identified with $\b^{q/2} \in D_q$ and
$\a^{p/2} \in D_p$ is identified with $\mu \in D_q$.}

\nd Remark:  In case (iv), $G(p,q)$ cannot, in general, be written
as an amalgamated free product of $\Z_p$ (or $D_p$ or $G(p,4)$)
with $\Z_q$ (or $D_q$ or $G(4,q)$).  There are simply too many
relations.  However, in this case $G(p,q)$ is equal to $G(m,4)$,
where $m=lcm(p,q)$.  This does turn out to be an amalgamated free product
$$ D_m  *\hs-.045 {\lower.5ex\hbox{}}_{D_4} G(4,4), \eqno(1.7) 
$$
where the $D_4$ subgroup is generated by $R_x^{\pi/2}$ and $R_z^\pi$.
In Theorem 2 we construct
a canonical form for $G(m,4,1)$, which can then be applied to this case.
\vs.1
\nd {\it Proof of the theorem:}
\  In cases (i), (ii) and (iii), we consider 
the natural map
$\rho$ from the abstract group $< \a, \b :\hbox{ (relations)}>$ to $G(p,q)$,
which sends $\a$ to $A$ and $\b$ to $B$.  In each case the relations
get mapped to the identity matrix $\I$, so the maps are well-defined.  Since
$A$ and $B$ are in the image of $\rho$, $\rho$ is onto.  We must
show that $\rho$ is 1--1.  To do this we show that every element $g \ne e$
of the abstract group can be written as a word in $\a$ and $\b$ in a
canonical way, and that the corresponding word in $A$ and $B$ is not
equal to the identity matrix.  The key to understanding these
words is the following lemma, which we assume at this point, and prove
after giving the rest of the proof of our theorem.

\nd {\bf Lemma 1:}\ {\it Let $m=s2^t$, $s$ odd and $t\ge 0$, and define 
$T=R_x^{2\pi/m}$.
If $W,E\in G(4,4,1)$, $4a_j\ne0$ (mod $m$), $b_j$  odd
and $n>0$, then} 
$$WS^{b_1}T^{a_1}S^{b_2}T^{a_2}\cdots S^{b_n}T^{a_n}E\ne \I.
\eqno (1.8)$$
\vs.1
We use Lemma 1  with $m=pq$. 
Note that $A=T^{q}$ and $B=S^3T^{p}S$.  An arbitrary
element of the abstract group is a word in $\a$ and $\b$.  Using
the relations $\a^p$ and $\b^q$ such an element $g$ can always be put in 
the form
$$\a^{\tilde a_1} \b^{\tilde b_2} \cdots \a^{\tilde a_n} \b^{\tilde b_n}, 
\eqno (1.9)$$
with each $\tilde a_i \in (0,p)$ and each $\tilde b_i \in (0,q)$, 
except $\tilde a_1$ and
$\tilde b_n$, which may equal zero.  In cases (ii) and (iii) we will use the
given relations to put further restrictions on the $\tilde a_i$'s and 
$\tilde b_i$'s.  The matrix $\rho(g)$
then equals
$$T^{a_1}S^3T^{b_1}S \cdots T^{a_n}S^3T^{b_n}S, \eqno (1.10)
$$
where $a_i = q \tilde a_i$ and $b_j = p \tilde b_j$.
We will show that, if $g \ne e$, this matrix cannot equal the identity.

i) If $p\ge 3$ and $q\ge 3$ are odd, we have $4 a_i\ne 0$ (mod $m$),
for $i>1$ and $4 b_j\ne 0$ (mod $m$) for $j<n$.  By Lemma 1, the only
way (1.10) can equal the identity is if $n=1$ and $a_1=b_1=0$, in which
case $g=e$ to begin with.  So $\rho$ is 1-1.

ii) Assume $p\ge 4$ is even and $q\ge 3$ is odd. By applying the
identity $\b^{b} \a^{p/2} = \a^{p/2} \b^{-b}$ to the expression (1.9)
we can require that all the $\tilde a_i$'s, except possibly $\tilde
a_1$, be nonzero and lie in the interval $(-{p \over 4},{p \over 4}]$.
Since $q$ is odd, we already have $4b_j\ne 0$ (mod $m$) for $j<n$.

Lemma 1 cannot be directly applied to expression (1.10), since, if $p$
is divisible by 4, $a_i$ may equal $m/4$ for some $i$.  We remove the
offending $T^{m/4}$ terms using the identity $S T^{m/4} S^3 = 
T^{m/4} S^3 T^{-m/4}$.  The factors of $T^{\pm m/4}$ may then be attached
to the neighboring $T^{b_j}$'s, producing new $b_j$'s that still satisfy
$4 b_j \ne 0 \pmod{m}$.  We then apply Lemma 1 to this expression.
As in case (i), the only way (1.10) can
equal the identity is if $n=1$ and $b_1=0$, in which case $\rho(g)=A^{a_1}$.
This equals the identity only if $a_1=0$, in which case $g=e$. So $\rho$ 
is 1-1.

iii) Assume $p\ge 4$ is even and $q=2s$, $s\ge 3$ odd.  By applying
the identity $\b^{q/2} \a^{a} =\a^{-a} \b^{q/2}$ we can require that
all the $\tilde b_i$'s, except possibly $\tilde b_n$, be nonzero and
lie in the interval $(-{q \over 4},{q \over 4})$.  By applying the
identity $\b^{b} \a^{p/2} = \a^{p/2} \b^{-b}$ to the expression (1.9)
we can require that all the $\tilde a_i$'s, except possibly $\tilde
a_1$, be nonzero and lie in the interval $(-{p \over 4},{p \over 4}]$.
As in the last case, we have $4 b_j \ne 0 \pmod{m}$ but may have $4
a_i = m$. The $T^{m/4}$ terms are eliminated as in case (ii), and
Lemma 1 shows that the only way $\rho(g)$ can equal zero is if $n=1$,
$4 \tilde a_1 = 0 \pmod{p}$ and $2 \tilde b_1 = 0
\pmod{q}$.  These 8 cases are easily listed, and the only one that gives
$\rho(g)= \I$ is $\tilde a_1 = \tilde b_1 = 0$, in other words $g=e$. 

iv) Since ${lcm}(p,q)$ is a multiple of $p$ and of $q$,
$G(lcm(p,q),4,1)$ contains $A$ and $R_x^{2\pi/q}$. Since it also
contains $S$, it contains $S^{-1}R_x^{2\pi/q}S=B$. Therefore
$G(lcm(p,q),4,1)$ contains $G(p,1,q)$.

Since 4 divides both $p$ and $q$, $G(p,1,q)$ contains $R_x^{2\pi/4}$
and $R_z^{2\pi/4}$, and so contains
$S=R_z^{-2\pi/4}R_x^{2\pi/4}R_z^{2\pi/4}$. But then it also contains
$SBS^{-1}=R_x^{2\pi/q}$, and so contains
$R_x^{2\pi/{lcm}(p,q)}$, and so contains
$G(lcm(p,q),4,1)$. (We have used the fact that
${lcm}(p,q)\,{ gcd}(p,q)=pq$, so there exist integers $k,\ell$
such that ${1\over{lcm}(p,q)}={k\over p}+ {\ell\over q}$). \qed

\vs.1
\nd {\it Proof of Lemma 1}:\ The lemma for a fixed value of $m$ is a corollary
of the lemma applied to $4m$.  So, without loss of generality, we may
assume from the start that $m$ is divisible by 4.  Let $x=e^{{2\pi
i/m}},\ y=x^s$ and $z=x^{2^t}$. Note that $y^{2^t}=1=z^s$. Since $s$
and $2^t$ are relatively prime in $\Z$, $\Z_m = \Z_s \times \Z_{2^t}$;
for each exponent $a$ there are $u,v\in
\Z$ such that $x^a= y^uz^v$. Let $R$ be the ring $\Z[x]=\Z[y,z]$. By
using the identity $y^{2^{t-1}}=-1$ we can write any element of $R$ in
the form \hfil 
$$\sum_{j=0}^{2^{t-1}-1}k_j(z)y^j, \hbox{ with }k_j(z)\in
\Z[z]. \eqno (1.11) 
$$  

\nd To see that this form is unique, we recall some facts about
the Euler function and cyclotomic polynomials.  The Euler function $\phi(n)$ 
gives the number of positive integers 
$r \le n$ relatively prime to $n$.  The cyclotomic polynomial of
$e^{2 \pi i/n}$ has order $\phi(n)$, so that $\Z[e^{2 \pi i/n}]$ has exactly
$\phi(n)$ linearly independent elements over $\Z$.  Now $\phi(s2^t)
= 2^{t-1} \phi(s)$, since for a number to be relatively prime to $s2^t$ it
must be odd and relatively prime to $s$. There are $\phi(s)$ such numbers
between 1 and $2s$, another $\phi(s)$ between $2s+1$ and $4s$, and so on.
But the form (1.11) requires exactly $2^{t-1}\phi(s)$ coefficients, $\phi(s)$
for each power of $y$.  If any of these could be eliminated, $\Z[x]$ would
be generated, as an abelian group, by fewer than $\phi(m)$ elements, which
is a contradiction.
\vs.1
Consider each factor $S^bT^a$ in the statement of the lemma. It is
of the form
$$ST^a=\pmatrix{0&-\tilde s&\tilde c\cr 0&\tilde c&\tilde s\cr -1&0&0\cr},
\eqno (1.12)$$
or
$$S^3T^a=\pmatrix{0&\tilde s&-\tilde c\cr 0&\tilde c&\tilde s\cr 1&0&0\cr},
\eqno (1.13)$$
where $\tilde c=\cos(2\pi a/m)=(x^a+\bar x^a)/2,\ \tilde s=\sin(2\pi a/m)
=(x^a-\bar x^a)/2i$. Writing $x^a$ in the form $x^a=y^uz^v$, we
distinguish each factor by whether $v=0 \! \pmod s$ or $v\ne 0 \! \pmod s$.
Let $I$ be a maximal extension of the ideal $(1+y) \subset R$.
\vs.1
We need two further lemmas, whose proofs we again defer.

\nd {\bf Lemma 2:} {\it If $v\ne 0 \! \pmod s$ when writing $x^a=y^uz^v$, 
the $(1,2)$, $(1,3)$, $(2,2)$ and $(2,3)$ entries of the matrix
$2S^bT^a$ are in $R$ but not in the maximal ideal $I$.}
\vs.1
\nd {\bf Lemma 3:} {\it If $x^a=y^{u}$, there is a power 
$w$ such that the $(2,2)$ entry of
$(1+y)^wS^bT^a$, namely $(1+y)^w(y^{u}+y^{-u})/2$,  
is in $R$ but not in the maximal ideal $I$.  
In particular, if $u= r2^k$, with
$r$ odd, then $w= 2^{t-1} - 2^{k+1}$.  Similarly, 
the $(1,2)$, $(1,3)$ and $(2,3)$ entries of $(1+y)^w S^b T^a$ are also
in $R$ but not in $I$.}
\vs.1
Now consider the matrix
$ F_i S^{b_i} T^{a_i}$,
where $F_i$ is either 2 or an appropriate power 
of $(1+y)$. (Note that $2\in I$ since $2=1-y^{2^{t-1}}=
(1+y)(1-y+y^2-\cdots -y^{2^{t-1}-1})$.)
We have shown that, modulo $I$, this matrix takes the form
$$ \pmatrix{0&\alpha&\beta\cr 0&\gamma&\delta\cr 0&0&0\cr},\eqno (1.14)$$
with $\alpha, \beta, \gamma, \delta$ nonzero elements of the field
$R/I$.  But the product of two (or more) matrices of this form again takes
this form, so   
$$FS^{b_1}T^{a_1}S^{b_2}T^{a_2}\cdots S^{b_n}T^{a_n},\eqno (1.15)$$ 
where $F$ is the appropriate product of the $F_i$'s, again takes this
form.  Matrices in the group $G(4,4,1)$ are, up to sign, permutation
matrices, so
$$FWS^{b_1}T^{a_1}S^{b_2}T^{a_2}\cdots S^{b_n}T^{a_n}E \eqno (1.16)
$$
has 4 matrix elements that are nonzero in $R/I$.  But $F$ times the
identity matrix is clearly zero modulo $I$, so
$WS^{b_1}T^{a_1}S^{b_2}T^{a_2}\cdots S^{b_n}T^{a_n}E$ can never equal
the identity.  \qed
\vs.1

\nd {\it Proof of Lemma 2}:  
We prove this first for the $(2,2)$ entry $x^a + x^{-a} =
y^uz^v+y^{-u}z^{-v}$.  Assume $y^uz^v+y^{-u}z^{-v}\in I$.  Since $1+y
\in I$, $(-y)^u-1 =-[1+y][1+(-y)+(-y)^2+\cdots +(-y)^{u-1}]
\in I$ and so $(-y)^u z^v - z^v \in I$.  Similarly, 
$(-y)^{-u}-1 \in I$,  so $(-y)^{-u} z^{-v} - z^{-v} \in I$.
This implies, using $y^uz^v+y^{-u}z^{-v}\in I$, that $z^v+z^{-v}\in I$. 
We now show that this implies $1\in I$, which is a
contradiction which proves the lemma for the $(2,2)$ entry. 

Let $\tilde z\equiv z^v \ne 1$.  Note that $\tilde z^s=1$. 
Now $(\tilde z+\tilde z^{-1})(\tilde z^2+\tilde z^3)=(\tilde z+\tilde
z^2+\tilde z^3+\tilde z^4)\in I$. Multiplying by $1+\tilde z^4+\tilde
z^8+\cdots +\tilde z^{4k}$ we see that $\tilde z+\tilde z^2+\tilde
z^3+\cdots +\tilde z^{4k+4}\in I$. We now consider two cases. If $s= 1$
(mod 4), take $k=(s-5)/4$, obtaining that $\tilde z+\tilde z^2+\tilde
z^3+\cdots +\tilde z^{s-1}\in I$. But $1+\tilde z+\tilde z^2+\tilde
z^3+\cdots +\tilde z^{s-1}=(1-\tilde z^{s})/(1-\tilde z)=0$, so $\tilde
z+\tilde z^2+\tilde z^3+\cdots +\tilde z^{s-1}=-1$, which implies
$1\in I$. Alternatively, if $s=3$ (mod 4) take $k=(s-3)/4$, obtaining
$\tilde z+\tilde z^2+\tilde z^3+\cdots +\tilde z^{s+1}\in I$. But
using $1+\tilde z+\tilde z^2+\tilde
z^3+\cdots +\tilde z^{s-1}=(1-\tilde z^{s})/(1-\tilde z)=0$,\
$\tilde z+\tilde z^2+\tilde z^3+\cdots +\tilde z^{s+1}=\tilde
z^{s+1}=\tilde z$, and if $\tilde z\in I$ then $1\in I$. 

Now consider the other entries.  The $(1,3)$ entry is just plus or minus 
the $(2,2)$ entry.  The $(1,2)$ and $(2,3)$ entries are (up to sign) of 
the form $y^{u'}z^v+y^{-u'}z^{-v}$, where $u' = u + 2^{t-2}$.  
The above argument, with $u$ replaced by $u'$, shows that these 
elements are in $R$ but not in $I$. \qed

\nd {\it Proof of Lemma 3}:  We essentially have to count the number $w$ of 
factors $(1+y)$ it takes 
so that $(1+y)^w(y^u + y^{-u})$ is a multiple of 2 in $\Z[y]$.
(It is important to note that this takes place in $\Z[y]$ not $\Z[x]$,
as we shall see.)
We first establish a few simple facts about powers of $(1+y)$. 
\item 1 If $c$ is a power of 2, then $(1+y)^c = 1 + y^c = 1 - y^c$
(mod 2). In particular, $(1+y)^{2^{t-1}} = 0 \pmod{2}$.  
(This follows from the binomial theorem).
\item 2 If $c$ is a power of 2, then $(1 \pm y^c)(1+y)^{2^{t-1}-c} = 0 
\pmod{2}$.  (This follows from 1, applied first to $c$ and then to
$2^{t-1}$.)
\item 3 If $c$ is a power of 2, then $(1 \pm y^c)(1+y)^{2^{t-1}-c-1} \not
= 0 \pmod{2}$. (The coefficient of $y^{2^{t-1}-1}$ is $\pm 1$, not a 
multiple of 2.)

Now we write
$$ \eqalign{ y^u + y^{-u} & = 
y^{-u}(1 + y^{2u}) \cr
& = y^{-u} (1 - y^{2^{k+1}} + y^{2^{k+1}} + y^{2u}) \cr
& = y^{-u} (1 - y^{2^{k+1}}) + y^{2^{k+1}-u}( 1 - y^{2u- 2^{k+1}}) + 2 y^u.}
\eqno (1.17)$$
The last term on the last line is always a multiple of 2.  Now 
$2u - 2^{k+1}={(r-1)\over 2}2^{k+2}$ so
$1 - y^{2u- 2^{k+1}}=
(1 - y^{2^{k+2}})(1+y+\cdots+y^{[{r-1\over 2}-1]2^{k+2}})$.
Therefore whenever $w \ge 2^{t-1} - 2^{k+2}$,
$(1+y)^w$ times the second term is divisible by 2.  But $(1+y)^w$ times 
the first term is divisible by 2 if and only if $w \ge 2^{t-1} - 2^{k+1}$.
As a result, 
$(1+y)^w (y^u + y^{-u})$ is divisible by 2 when $w = 2^{t-1} - 2^{k+1}$,
but is not divisible by 2 when $w=2^{t-1} - 2^{k+1}-1$.

Now let  $u' = u + 2^{t-2}$, as before.  Since $k$ is, by assumption, 
less than $t-2$, the power of $(1+y)$ needed
to make $y^{u'} + y^{-u'}$ divisible by 2 is the same as that
needed to make $y^u + y^{-u}$ divisible by 2.

So we have determined the critical power $w$ such that multiplying by
$(1+y)^w$ puts the matrix elements in $\Z[y]$ but not in the ideal
$(1+y)_0$ in $\Z[y]$ generated by $1+y$. $(1+y)_0$ is a maximal ideal
in $\Z[y]$ since $\Z[y]/(1+y)_0$ is the field $\Z_2$.
Since $I\cap\Z[y]$ must be a proper ideal in $\Z[y]$,
$I\cap\Z[y]$ must coincide with $(1+y)_0$, and so the matrix elements
cannot be in $I$.
\qed

This completes the proof of Theorem 1.
\vs.3
\nd {\bf \S 2.\ Canonical forms for G(p,q)}
\vs.1
In this section we construct canonical forms for elements of the
groups $G(p,q)$.  Since $G(p,q)$ is always a subgroup of $G(pq,4,1)$, we
first construct a canonical form for elements of $G(m,4,1)$, 
where $m$ is an arbitrary
integer (Theorem 2).  This is most useful when
$p$ and $q$ are both divisible by 4, for in that case 
$G(p,q)=G(lcm(p,q),4,1)$.
In the remaining cases, where $p$ or $q$ is not divisible by 4, Theorem~3
provides canonical forms for
elements of $G(p,q)$ as products of the generators of $G(p,q)$.
 
As before, we take $S=R_y^{2\pi/4}$ and 
$T=R_x^{2\pi/m}$.  We also define $U=R_x^{2\pi/4}$.  Note that $S$
and $U$ generate $G(4,4,1)=G(4,4,4)$.
\vs.1
\nd {\bf Theorem 2: Canonical form for G(m,4,1)}

{\it Let $H= G(4,4,1) \cap G(m,4,1)$.  Let $g$ be an arbitrary element
of $G(m,4,1)$.  Then $g$ can be uniquely written
in the form
$$ g = W ST^{a_1} \cdots ST^{a_n} E,  \eqno(2.1) $$
for some  $n \ge 0$, where $W$ and $E$ are elements of $H$, 
$a_i$ is an integer, and the following restrictions are applied:
\item {1.} If $m$ is odd, $W \in \{\I, S^3 \}$, $a_i \in (-{m \over 2},
{m \over 2})$, and $a_i \ne 0$.
\item {2.} If $m$ is twice an odd number, $W \in \{\I, S^3 \}$, 
$a_i \in (-{m \over 4},{m \over 4})$ and $a_i \ne 0$.
\item {3.} If $m$ is divisible by $4$, $W \in \{\I, S^3, U \}$, 
$a_i \in (-{m \over 4},{m \over 4})$, $a_i \ne 0$, and $a_n \in 
(0,{m \over 4})$.
\item {4.} If $n=0$, then $W=\I$.
}

\nd Remark 1:  Since $G(m,4,1) \subset G(2m,4,1) \subset G(4m,4,1)$, one
could write any element $g$ of $G(m,4,1)$ using the canonical form for 
$G(4m,4,1)$. However, if $m$ is not divisible by 4, this would typically 
involve writing $g$ as a product of matrices that are not themselves
in $G(m,4,1)$.

\nd Remark 2:  The allowed values of $W$ can be understood as follows.  
There is a subgroup $H_1$ of $H$ that can be commuted (or anticommuted)
past powers of $ST^a$, or absorbed into $ST^a$.  
Factors in $H_1$ can be removed from $W$ and either absorbed into 
$ST^{a_1}$ or transferred all the way from $W$ to $E$.  
The allowed values of $W$ are representatives of the cosets in $H/H_1$.

If $m$ is odd, then $H = \Z_4$, generated by $S$.  In this case
$H_1=\Z_2$, generated by $S^2$, and $H = H_1 \cup S^3 H_1$.  If $m$
is twice an odd integer, then $H$ is the 8-element group generated by
$S$ and $U^2$, $H_1$ is the 4-element subgroup generated 
by $S^2$ and $U^2$, and once again $H = H_1 \cup S^3 H_1$.  
If $m$ is divisible by 4, then
$H=G(4,4,1)$ is the 24-element group generated by $S$ and $U$, 
$H_1$ is the 8-element
subgroup generated by $S^2$ and $SUS^{-1}$, and $H= H_1 \cup S^3 H_1
\cup U H_1$.

\nd Remark 3: When $m$ is divisible by 4, the canonical form (2.1) is
closely related to the amalgamated free product (1.7).  The nontrivial
cosets of $G(4,4)/D_4$ are represented by $S$ and
$SU^{-1}$, while the nontrivial cosets of $D_m/D_4$ are represented by $T^a$, 
with $a \in (0,{m\over 4})$.  Multiplying these together we get
$ST^a$, with $a$ nonzero and in $(-{m\over 4},{m\over 4})$.

\nd {\it Proof}: The proof is an application of two lemmas, which are
proved below:

\nd {\bf Lemma 4}:  {\it Any element of $G(m,4,1)$ can be put in the form 
$(2.1)$.}  

\nd {\bf Lemma 5}:  {\it If $g \in G(m,4,1)$ is in the form $(2.1)$, 
there is at most one expression $g'$ in the form $(2.1)$ such that
that $g g'=\I$.}

By Lemma 4, representatives for $g$ and $g^{-1}$ always exist.
Applying Lemma 5 to $g^{-1}$ we see that the representation for $g=
(g^{-1})^{-1}$ is unique, and the theorem is proved. \qed

\nd {\it Proof of Lemma 4}:  There are 3 cases to consider, depending on 
whether $m$ is odd, twice an odd number, or divisible by 4.  In all cases
we assume that $g$ is not in $H$, since if $g \in H$ we can simply take
$W=\I$, $n=0$, $E=g$.

Let $m$ be odd.  Any element $g$ of $G(m,4,1)$ can be written as a
word in the generators $S$ and $T$, and hence takes the form
$$ g = S^{b_1} T^{a_1} S^{b_2} T^{a_2} \cdots S^{b_N} T^{a_N} S^{b_{N+1}}, 
\eqno(2.2) $$
with no restrictions on $N$ or $b_i$ or $a_i$.  By applying the
relations $S^4=T^m=\I$, we can force each $a_i \in (-{m \over 2},{m
\over 2})$, $b_i \in \{0,1,2,3\}$.  If any $a_i=0$ or $b_i=0$, we can
collapse the expression into a shorter word and proceed as before.  If
any $b_i=2$, we can use the relation
$$ S^2 T^{a} = T^{-a} S^2 \eqno(2.3) 
$$ 
to shorten the word further.  
Since the word has finite length, this process must terminate, leaving us
with an expression of the form
$$ g = S^{b_1} T^{a_1} \cdots S^{b_n} T^{a_n} S^{b_{n+1}}, \eqno (2.4) 
$$
where each $b_i$, with the possible exception of $b_1$ and $b_{n+1}$,
is odd, and each $a_i \in (-{m \over 2},{m \over 2})$ and is nonzero.
Next we force $b_1$ to equal $0$ or $1$ by using (2.3), if needed, to
push a factor of $S^2$ past $T^{a_1}$.  We then force $b_2=1$ by
possibly using (2.3) to push a factor of $S^2$ past $T^{a_2}$.
Continuing in this way we can make all the $b_i$'s equal to $1$, with
the possible exceptions of $b_1$, which can equal $0$ or $1$, and
$b_{n+1}$, which is not constrained.  Define $W=S^{b_1-1}$,
$E=S^{b_{n+1}}$.  Our element $g$ then takes the form (2.1).

Note that the specific numbers $a_i$ may be changed in converting from
the form (2.2) to (2.4) to (2.1).  In our usage $a_i$ does not denote
a fixed number; rather, it denotes the $i$-th exponent of $T$ in a
typical expression.

Now suppose that $m$ is twice an odd number.  We proceed as before to
reach the form (2.4), with $b_i$ odd and $a_i \in (-{m \over 2},{m
\over 2}]$ and nonzero.  We then use the relation
$$T^{m/2} S^b = S^{4-b} T^{m/2}, \eqno (2.5) 
$$ 
as needed, to eliminate factors of $T^{m/2}$ and to make each $a_i \in
(-{m \over 4},{m \over 4})$.  We begin with $a_1$, possibly using
(2.5) to push $T^{m/2}$ past $S^{b_2}$, then forcing $a_2$ into $(-{m
\over 4},{m \over 4})$ by possibly 
using (2.5) to push $T^{m/2}$ past $S^{b_3}$, and so on.  In this way
all the $a_i$'s, with the possible exception of $a_n$, can be put in
$(-{m \over 4},{m \over 4})$.  Then we use (2.3) to make all the
$b_i$'s, with the same exceptions for $b_1$ and $b_{n+1}$ as before,
equal to 1.  Note that minus a nonzero integer in $(-{m \over 4},{m
\over 4})$ is another nonzero integer in $(-{m \over 4},{m \over 4})$,
so fixing the $b_i$'s does not disrupt the form of the $a_i$'s.

This gives us an expression of the form (2.4), with each $a_i \in (-{m
\over 4},{m \over 4})$ and nonzero, except $a_n$, which is nonzero and
in $(-{m \over 2},{m \over 2})$, and with each $b_i=1$, except for
$b_1$ which may equal 0 or 1, and $b_{n+1}$ which is arbitrary. As
before, define $W=S^{b_1-1}$.  If $a_n \in (-{m \over 4},{m \over
4})$, define $E=S^{b_{n+1}}$; otherwise, define $E = T^{m/2}
S^{b_{n+1}}$.  This puts us in the form (2.1).

Finally, suppose that $m$ is divisible by $4$.  We proceed as before
to the form (2.4), with $b_i$ odd and $a_i$ not divisible by $m/2$.
If any of the $a_i$'s (other than $a_1$ or $a_n$) is divisible by
$m/4$, we can reduce the length of the word further, as follows.
First use (2.3) to set $b_{i}=b_{i+1}=1$.  Then use one of the
relations
$$ S U S = U S U; \qquad \qquad S U^3 S = U^3 S U^3 \eqno (2.6) 
$$ 
to change $T^{a_{i-1}} S T^{\pm m/4} S T^{a_{i+1}}$ to $T^{a_{i-1} \pm
{m \over 4}} S T^{a_{i+1} \pm {m \over 4}}$.  This may result in an
exponent that is divisible by $m/2$, in which case we use (2.5) to
reduce the word length further.  Since the original word has finite
length, we eventually reach the form (2.4) where none of the $a_i$'s,
possibly excepting $a_1$ and $a_n$ is divisible by $m/4$.  

If $a_1$ is divisible by $m/4$ we define $\tilde W=S^{b_1} T^{a_1}$;
otherwise $\tilde W = S^{b_1-1}$.  If $a_n$ is divisible by $m/4$ we
define $\tilde E = S T^{a_n} S^{b_{n+1}}$; otherwise $\tilde E =
S^{b_{n+1}}$.  In any case, we now have $g$ in the form
$$\tilde W S^{b_1} T^{a_1} \cdots S^{b_n} T^{a_n} \tilde E, \eqno(2.7) 
$$
with $\tilde W$ and $\tilde E$ in $H=G(4,4,4)$, with $b_i$ odd and
with $a_i$ not divisible by $m/4$.

This is almost of the form (2.1). To achieve the necessary restrictions on
$W$, $E$, $a_i$ and $b_i$, we work from left to right, pushing undesired
factors rightwards.   Let $H_1$ be the 8-element
subgroup of $H$ generated by $S^2$ and $S U S^{-1}$.  Of these two
generators, $S^2$ can be commuted past a factor $S^bT^{a}$ (changing it to
$S^b T^{-a}$), while $SUS^{-1}$ can be absorbed into a factor of $ST^{a}$:
$$
SUS^{-1}S T^{a} = SU T^{a} = ST^{a+{m\over 4}}, \eqno(2.8)
$$
without changing the fact that $a$ is not divisible by ${m\over 4}$.
Factors in $H_1$ can thus be
removed from $\tilde W$ and moved rightwards.  Since $H = H_1 \cup S^3
H_1 \cup U H_1$, we can change $\tilde W$ to $\I$, $S^3$ or $U$, which
we then call $W$.  Then, working left to right, we use (2.3) to change
some $b_i$'s from 3 to 1 and use $(2.5)$ to place the $a_i$'s in the
range $(-{m\over 4},{m\over 4})$.  Finally, if $a_n<0$, we define $E$
to be $U^{-1} \tilde E$ (otherwise $E=\tilde E$).  By factoring out
$U^{-1}$, we put $a_n \in (0,{m \over 4})$, and we have achieved the
form (2.1).  \qed

\nd {\it Proof of Lemma 5}: Suppose that we have 
$g = W S T^{a_1} \cdots S T^{a_n} E$, and that $g^{-1} = W' S T^{a_1'}
\cdots S T^{a'_{n'}} E'$, with appropriate restrictions on $W, a_i, E,
W', a'_i, E'$.  We will show that there is a unique choice of $W'$,
$n'$, $a_i'$ and $E'$.  Any other choices will allow us to turn the
expression $W S T^{a_1} \cdots S T^{a_n} E W' S T^{a_1'} \cdots S
T^{a'_{n'}} E'$ into something of the general form
$$ W S^{b_1} T^{a_1} \cdots S^{b_N} T^{a_N} E, \eqno(2.9) $$ 
with $b_i$ odd, $a_i$ not divisible by $m/4$, and $W,E \in H$.  By
Lemma 1, such an expression is not equal to $\I$, contradicting the
equation $g g^{-1} = \I$.  As usual, the details depend on whether $m$
is odd, twice odd, or divisible by 4.

Suppose $m$ is odd.  We must choose $W'$ such that $E W' S$ is an even
power of $S$. If this choice is not made, then $W S T^{a_1} \cdots S
T^{a_n} (E W' S) T^{a'_1} \cdots S T^{a'_{n'}} E'$ is of the form
(2.9).  Since $W' \in \{\I, S^3 \}$, there is exactly one right
choice.

Now, since $E W' S$ is an even power of $S$, it can be commuted past
all the $S T^{a'}$ factors, leaving us with the form $W S T^{a_1}
\cdots S T^{a_n} T^{a'_1} \cdots S T^{a'_{n'}} E'$, where the new
$a'$s are $\pm$ the old ones, and the new $E'$ is $EW'S$ times the old
one.  If $a_n + a'_1 \ne 0$, then we are again of the form (2.9), so
we must have $a'_1 = -a_n$.  We again push a factor of $S^2$ all the
way to the right, and find that there is a unique value of $a'_2$ such
that we again avoid the form (2.9).  This process continues, with each
$a'_i$ determined by $a_{n+1-i}$ and the history of what has passed
before.  We cannot have $n' \ne n$, as that would leave some powers of
$ST^a$ (or $ST^{a'}$) that are not cancelled.  If $n = n'$ and each
$a_i'$ is chosen correctly, we eventually reach the form $W \times
$(transferred powers of $S) \times E'$.  There is clearly a unique
choice of $E'$ that makes this equal unity.

Now suppose $m$ is twice an odd number.  The argument is almost
identical.  $E$ is a power of $S$, possibly times $U^2$.  As before,
if that power is odd, we must choose $W'=\I$, while if that power is
even we must choose $W'=S^3$.  If this choice is not made, we can
commute any $U^2$ factors to the right and achieve the form (2.9),
which would not be the identity.  If this choice {\it is} made, then
$E W' S$ is an even power of $S$, possibly times $U^2$, and can be
commuted past all the $S T^{a'}$ factors.  The argument then proceeds
precisely as before, with each $a_i'$ determined by $a_{n-i+1}$, and
with $E'$ determined by what remains after the $ST^{a}$ factors are
all cancelled.

The same line of reasoning works for $m$ divisible by 4, with a few
extra steps to deal with complications coming from powers of $U$.
Recall that we have the 8-element subgroup $H_1$ of $H$, generated by
$SU S^{-1}$ and $S^2$, of elements that can be commuted past (or
absorbed into) $S T^a$.  We write $g g^{-1}$ as $W S T^{a_1} \cdots S
T^{a_n} S^{-1} (S E W') S T^{a'_1}
\cdots S T^{a'_{n'}} E'$.  $S E W'$ can be expressed as $xh$, where $x
\in \{ \I, S^3, U \}$ and $h$ is an element of $H_1$.  We can push $h$
all the way to the right, getting an expression of the form
$$W  S T^{a_1} \cdots S T^{a_n} S^{-1} x S T^{a'_1} \cdots S T^{a'_{n'}} E'.
\eqno (2.10) $$
If $x=S^3$, this is of the form (2.9) and cannot equal the identity.
If $x=U$, we use the identity $S^{-1} U S = U^3 S^3 U$ and absorb the
powers of $U$ into $T^{a_n}$ and $T^{a_1'}$ to put this in the form
(2.9).  Thus the only way to have $g g^{-1} =\I$ is to have $x=\I$, or
equivalently for $S E W' \in H_1$.  It is straightforward to check
that, for each possible $E \in H$, there is a unique $W' \in \{ \I,
S^3, U \}$ such that $S E W' \in H_1$.

Once $W'$ is chosen and $h$ is pushed to the right, we have an
expression of the form $W S T^{a_1} \cdots S T^{a_n} T^{a'_1} \cdots S
T^{a'_{n'}} E'$, where the new $a'$s and $E'$s are determined in a
1--1 way by the old ones.  If $a_n + a_1'$ is not divisible by $m/4$,
this is of the form (2.9) and cannot equal unity.  Since $a_1' \in
(-{m \over 4},{m \over 4})$, there are exactly two values of $a_1'$
for which $a_n + a_1'$ is divisible by $m/4$, one of which has $a_n +
a_1' = 0$, the other of which has $a_n + a_1'= \pm m/4$.  If $a_n +
a_1'= \pm m/4$, we use the identity (2.6), and the fact that neither
$a_{n-1}$ nor $a_2'$ is a multiple of $m/4$, to achieve the form
(2.9).  Thus we must have $a_1' = -a_n$.

Similarly, $a_2'$ is determined by $a_{n-1}$, and so on.  As before,
we must have $n = n'$.  After $n-1$ cancellations we are left with
$$W S T^{a_1} T^{a_n'} \hbox{(transferred powers of $S^2$)} h E'. 
\eqno (2.11) $$
At this point the argument that $a_1 + a_n' \ne \pm m/4$ breaks down.
However, $a_n'$ is restricted to $(0,{m \over 4})$. Either $-a_1$ or
${m \over 4}-a_1$, but not both, lie in $(0,{m \over 4})$.  This is
the only possible value of $a_n'$ that keeps (2.11) from being of the
form (2.9).  Once this choice is made, $W S T^{a_1} T^{a_n'}
$(transferred powers of $S^2) h \in H$, and $E'$ must be chosen to be
the inverse of this element.  \qed

\vs.1
We now turn to canonical forms for $G(p,q)$ in general.  As before,
let $A = R_x^{2 \pi/p}$ and let $B = R_z^{2 \pi/q}$.  If $p$ and $q$
are both odd, then $G(p,q)$ is a free product, and every element can
be uniquely written in the form
$$ A^{a_1} B^{b_1} \cdots A^{a_n} B^{b_n}, \eqno(2.12) 
$$
with $a_i \in (-{p \over 2},{p \over 2})$, $b_j \in (-{q \over 2},{q
\over 2})$, and all exponents, except perhaps $a_1$ and $b_n$,
nonzero.  If $p$ and $q$ are both divisible by 4, then $G(p,q) =
G({lcm}(p,q),4,1)$, and a canonical form is provided by Theorem 2.
But what if $p$ is even and $q$ is not divisible by 4?  Here we define
three canonical forms in such cases.  Depending on the application,
one or another of these forms may be most useful.
\vs.1
\nd {\bf Definition:} 
{\it Let $p$ and $q$ be positive integers $\ge 3$ with $p$ even and
$q$ not divisible by $4$. 

A product $(2.12)$, with all exponents nonzero except perhaps $a_1$
and $b_n$, is in {L-canonical form} if: For $i>1$, $a_i \in (-{p \over
4},{p \over 4}]$; $a_1 \in (-{p \over 2},{p \over 2}]$; $b_1 \in (-{q
\over 2},{q \over 2}]$, and may equal $q/2$ only if $n=1$; for $j>1$,
$b_j \in (-{q \over 2},{q \over 2})$ if $q$ is odd and $b_j \in (-{q
\over 4},{q \over 4})$ if $q$ is even.

A product $(2.12)$, with all exponents nonzero except perhaps $a_1$
and $b_n$, is in {R-canonical form} if: For $i<n$, $a_i \in (-{p \over
4},{p \over 4}]$; $a_n \in (-{p \over 2},{p \over 2}]$, and may equal
$p/2$ only if $n=1$; $b_n \in (-{q \over 2},{q \over 2}]$; for $j<n$,
$b_j \in (-{q \over 2},{q \over 2})$ if $q$ is odd and $b_j \in (-{q
\over 4},{q \over 4})$ if $q$ is even.

A product $(2.12)$, with all exponents nonzero except perhaps $a_1$
and $b_n$, is in {C-canonical form} if: For $i>1$, $a_i \in (-{p \over
4},{p \over 4}]$; $a_1 \in (-{p \over 2},{p \over 2}]$; $b_n \in (-{q
\over 2},{q \over 2}]$; for $j<n$, $b_j \in (-{q \over 2},{q \over
2})$ if $q$ is odd and $b_j \in (-{q \over 4},{q \over 4})$ if $q$ is
even.}

The differences between the canonical forms is just a matter of where
to put factors of $R_x^\pi$ and $R_z^\pi$.  In L-canonical form they
are placed at the left, in R-canonical form they are placed at the
right, and in C-canonical form $R_x^\pi$ is moved left and $R_z^\pi$
is moved right.  If $q$ is odd, $R_z^\pi$ does not appear, and the L-
and C-canonical forms coincide.
\vs.1
\nd {\bf Theorem 3: Canonical forms for $G(p,q)$}

{\it Let $p$ and $q$ be positive integers $\ge 3$ with $p$ even and
$q$ not divisible by $4$.  Each element of $G(p,q)$ can be uniquely
written in L-canonical form, and uniquely written in R-canonical form,
and uniquely written in C-canonical form.}

\nd {\it Proof}:  The proof has several steps, 
and is quite similar in spirit to the proof of Theorem 2.  First we
show that any element of $G(p,q)$ can be put into each of the
canonical forms.  Next we show that the only R-canonical form for the
identity element is $A^0B^0$.  Then we show that R-canonical forms are
unique, by showing that each element in L-canonical form has a unique
inverse in R-canonical form.  Finally we show that L-canonical and
C-canonical forms are unique by relating them to R-canonical forms of
the same length.

\nd {\it Step 1:} Since $A$ and $B$ generate $G(p,q)$, and since
$A^p=B^q=\I$, any element of $G(p,q)$ can be written in the form
(2.12) with each $a_i \in (-{p \over 2}, {p \over 2}]$, each $b_j \in
(-{q \over 2},{q \over 2}]$, and with all terms except $a_1$ and $b_n$
nonzero.  If any $a_i$'s except $a_0$ equal $p/2$, we can shorten the
word using the identity $A^a B^b A^{p/2} B^{b'} = A^{a + {p \over 2}}
B^{b'-b}$.  Similarly, if any $b_j$ other than $b_n$ equals $q/2$, we
can shorten the word with the identity $A^a B^{{q \over 2}} A^{a'}
B^{b} = A^{a -a'} B^{b+{q \over 2}}$.  Thus we can achieve the form
(2.12) in which for $i>1, 0 \ne a_i \in (-{p \over 2},{p \over 2})$
and for $j <n, 0 \ne b_j \in (-{q \over 2},{q \over 2})$.

Suppose $q$ is odd.  To put our expression in L-canonical (or
C-canonical) form, we must adjust the exponents $a_j$, $a>1$ that are
not in $(-{p \over 4},{p \over 4}]$ by $\pm p/2$, using the identity
$A^{a_{i-1}} B^b A^{a_i \pm {p \over 2}} = A^{a_{i-1} + {p \over 2}}
B^{-b} A^{a^i}$.  We begin by adjusting $a_n$ at the expense of
$a_{n-1}$ and $b_{n-1}$, then adjust $a_{n-1}$ at the expense of
$a_{n-2}$ and $b_{n-2}$, and continue until all the $a_i$'s, except
perhaps $a_1$, are in $(-{p \over 4},{p \over 4}]$.  In the process
some of the $b_i$'s may change sign, but this does not change the
conditions $b_j \in (-{q \over 2},{q \over 2})$, $\b_j \ne 0$.  Also
note that, since $b_1$ was originally not divisible by $q/2$, it
remains not divisible by $q/2$.

To put an expression in R-canonical form one first adjusts $a_1$ at
the expense of $b_1$ and $a_2$, then adjusts $a_2$, and so on through
$a_{n-1}$.  It should be clear that an expression in L-canonical form
can be converted to an expression in R-canonical form {\it of the same
length}, and vice-versa.

Now suppose $q$ is even.  To put an expression in any of the canonical
forms, one first adjusts the exponents $b_j$ using the identity
$B^{b_i} A^a B^{b_{i+1} + {q \over 2}} = B^{b_i + {q \over 2}} A^{-a}
B^{b_{i+1}}$, so that all but the first (for L-canonical) or last (for
C-canonical or R-canonical) $b_i$ lie in $(-{q \over 4},{q \over 4})$.
One then adjusts the $a_i$'s, as above.  Since the condition $b_j \in
(-{q \over 4},{q \over 4})$ is equivalent to $-b_j \in (-{q \over
4},{q \over 4})$, adjusting the $a_i$'s does not disrupt the form of
the $b_j$'s.  Again, it should be clear that converting from one
canonical form to another does not change the length of the word.

\nd {\it Step 2:} We must show that a nontrivial word in R-canonical
form cannot equal the identity.  This is essentially a repeat of an
argument in the proof of Theorem 1.  Embed $G(p,q)$ in $G(m,4,1)$,
where $m=pq$.  $G(m,4,1)$ is generated by $T=R_x^{2\pi/m}$ and
$S=R_y^{\pi/2}$.  Note that $A=T^{q}$ and $B=S T^{p} S^3$.  Rewrite
any nontrivial word in $A$ and $B$ as a word in $S$ and $T$.  Although
the powers of $S$ in this word are all odd, the expression is not
quite in the form (2.9), as some $T^{m/4}$ factors may appear.  These
are removed with the identity $T^b S^3 T^{m/4} S T^{b'}= T^{b- {m
\over 4}} S^3 T^{b'+{m \over 4}}$.  Unless the original word was
$A^{(0 \hbox{\sevenrm{} or } {p\over 2})} B^{(0 \hbox{\sevenrm{} or }
{q \over 2})}$, the result is of the form (2.9), and by Lemma 1 is not
the identity.  The three special cases $A^{p/2}B^0$, $A^{0}B^{q/2}$
and $A^{p/2}B^{q/2}$ are separately checked to not equal the identity.

\nd {\it Step 3:} We show that every word in L-canonical form has
at most one inverse in R-canonical form.  We write $g = A^{a_1}
B^{b_1} \cdots A^{a_n} B^{b_n}$, $g^{-1}= A^{a_1'} B^{b_1'}
\cdots A^{a'_{n'}} B^{b'_{n'}}$, where $g$ is in L-canonical form
and $g^{-1}$ is in R-canonical form.  We show that, unless the $a'$s
and $b'$s are all chosen correctly, the product $g g^{-1}$ can be placed
in a nontrivial R-canonical form, and so cannot equal unity.  We
proceed by induction.

Uniqueness of inverses is easy to check for $n=1$.  The unique inverse
for $A^a B^b$, with $a$ and $b$ both nonzero, is $A^0 B^{-b} A^{-a}
B^0$, unless $q$ is even and $b \not \in (-{q \over 4},{q \over 4})$,
in which case the unique inverse is $A^0 B^{{q \over 2} -b} A^a B^{{q
\over 2}}$.  The unique inverse to $A^a B^0$ is $A^{-a} B^0$, the
unique inverse to $A^0 B^b$ is $A^0 B^{-b}$, and the unique inverse to
$A^0 B^0$ is $A^0 B^0$.

Now assume the assertion is proved for $n=k$ and that we have an
expression $g$ of length $n=k+1$.  If $b_n =0$, we must have $a_1' =
-a_{n}$ (unless $a_n=p/4$, in which case $a_1'=p/4$), as $g g^{-1} =
A^{a_1} B^{b_1} \cdots A^{a_{n}+ a'_1} B^{b'_1} \cdots A^{a'_{n'}}$.
If $a_1'$ is chosen incorrectly, the exponent $a_{n}+ a'_1$ is not
divisible by $p/2$, so, by transferring powers of $A^{p/2}$ and
$B^{q/2}$ from left to right, this expression can be placed in
nontrivial R-canonical form, and so cannot equal $\I$.  Once $a_1'$ is
chosen correctly, $g A^{a_1'}$ is a word of length $n-1=k$. It may be
converted to L-canonical form and so, by the inductive hypothesis, its
inverse in R-canonical form is uniquely determined.  But its inverse
is precisely $A^0 B^{b_1'} \cdots A^{a_{n'}'} B^{b_{n'}'}$, so $b_1$,
$a_2$, etc., are uniquely determined.

If $b_n \ne 0$, we must have $a_1 =0$, or else $g g^{-1} = A^{a_1}
B^{b_1} \cdots A^{a_{n}} B^{b_n} A^{a'_1} B^{b'_1} \cdots A^{a'_{n'}}$
could similarly be massaged into nontrivial R-canonical form.  By the
same argument, we also must have $b'_1 = - b_n$.  But then $g
B^{b'_1}$ is a word of length $k+1$ with final exponent zero, and its
inverse is uniquely determined by the argument of the previous
paragraph.

Since every element $g \in G(p,q)$ can be put in L-canonical form,
this shows that each element $g \in G(p,q)$ has a unique inverse in
R-canonical form.  Thus $g = (g^{-1})^{-1}$ has a unique R-canonical
form.

\nd {\it Step 4:} We count the number of words of length $n$ for any 
particular form by multiplying the number of choices for each $a_i$
and each $b_j$.  This number equals $p q \left ({p \over 2} - 1 \right
)^{n-1} (q-1)^{n-1}$ if $q$ is odd, and $p q \left ({p \over 2} -1
\right )^{n-1} \left ({q \over 2} -1 \right )^{n-1}$ if $q$ is even,
and is the same for R-canonical, L-canonical and C-canonical forms.

Now consider the set of R-canonical forms of length $n$ or less. We
have already shown that these can be converted to L-canonical forms of
length $n$ or less. Since the number of L-canonical forms equals the
number of R-canonical forms, and since each R-canonical form
corresponds to a distinct element of $G(p,q)$, each L-canonical form
is achieved in this way exactly once.  Thus distinct L-canonical forms
of length $n$ or less correspond to distinct R-canonical forms, and
hence to distinct elements of $G(p,q)$.  Since $n$ is arbitrary, this
shows that L-canonical forms are unique.  A similar argument shows
that C-canonical forms are unique.
\qed
\vs.3
\nd {\bf \S 3.\  Dite and kart tilings}
\vs.1
We construct here 3-dimensional tilings with symmetry group $G(10,4)$,
based on a version of the 2-dimensional kite and dart tilings.  The
new tilings consist of congruent copies of 8 elementary prisms,
constructed as follows.

Consider the two right triangles of Fig.\ 1 which we denote by
$\delta$ and $\kappa$.  $\delta$ has legs of lengths 1 and $\tau\sqrt
{2+\tau}$ and $\kappa$ has legs of lengths $\tau$ and $\sqrt
{2+\tau}/\tau$, where $\tau=(1+\sqrt 5)/2$, the golden mean. (The
triangles $\delta$ and $\kappa$ are halves of the triangles $S_A$ and
$S_B$ introduced by Raphael Robinson in his version of the kite and
dart tilings [GrS].) It is elementary to check that the small angle in
$\delta$ is $\pi/10$ and the small angle in $\kappa$ is $\pi/5$. We
next introduce triangles $\tilde \delta$ and $\tilde \kappa$ which are
larger than $\delta$ and $\kappa$ by a linear factor $\tau$.

The constructions in Fig.\ 1 (called deflation rules) then show how
these 4 triangles can be decomposed into congruent copies of triangles
which are each a linear factor $\tau^2$ smaller than $\delta,\ \kappa,\
\tilde \delta$ and $\tilde \kappa$.

We thicken $\delta$ by two different depths to make two types of prisms, 
a ``short thin dite'' of depth 1 and a ``short thick dite'' of depth
$\tau$. Likewise from $\tilde \delta$ we make a ``{\it tall} thin
dite'' of depth 1 and a ``{\it tall} thick dite'' of depth
$\tau$. Finally, replacing $\delta$ by $\kappa$ we make the analogous
4 types of ``karts''.

We now make deflation rules for the prisms as follows, 
again shrinking by a linear factor of
$\tau^2$.  We begin with the karts.

The short thin kart is deflated into a pair of layers. The ``top''
layer consists of short thin dites and karts which, when viewed from above,
have the same pattern as the deflated 2-dimensional $\kappa$ (Fig.\ 1). The
bottom layer consists of short thick
dites and karts in the same pattern. Since $\tau^2=1+\tau$, 
the sum of the thicknesses of the two layers equals the thickness of
the original thin kart.

The rule for the short {\it thick} kart is similar, only we now use 3
layers of short dites and karts; a top thin layer and 2
thick lower layers. Since $\tau=(1+2\tau)/\tau^2$, the total thickness 
of the deflated layers equals the thickness of the original thick kart.

The rules for the {\it tall} thin and thick karts are now immediate,
replacing the short dites and karts in the deflation of $\tilde
\kappa$ by tall dites and karts.

The rule for deflating dites uses the deflation of $\delta$ and $\tilde
\delta$ rather than that of $\kappa$ and $\tilde \kappa$, with one added
twist, based on the rectangles appearing in the deflations of $\delta$
and $\tilde \delta$ (Fig.\ 1).  As with karts, the thin dites deflate
into 2 layers and the thick ones into 3 layers.  The deflation of each
dite generates 2 or 3 parallelpipeds corresponding to the
aforementioned rectangle. If the original dite is short, then the
parallelpiped in the thin layer has 2 square faces, while if the
original dite is tall, the parallelpipeds in the thick layers have 2
square faces.  We then rotate the parallelpiped by $\pi/2$ about the
axis joining the centers of the square faces, as in Fig.\ 2.  This
completes the deflation rules.

Given these deflation rules, the dite and kart tilings are obtained as
follows. We begin with any one of the 8 prisms, say a short thick
dite, deflate it, then expand the 10 resulting small prisms linearly by
a factor $\tau^2$ about some point. Next reposition the 10 prisms so that
a copy of a short thick dite is sitting over the original
position. By indefinitely repeating this process of
deflation-expansion-repositioning one obtains the desired tilings of
space.

What is the group of relative orientations for each species of tile?
Since each species of tile, when deflated several times, gives rise
to all species of tiles, the group does not depend on the species.
We show that this group is $G(10,4)$. 

In the 2-dimensional $\delta$-$\kappa$ tiling, the group is $D_{10}$.
One can see two $\delta$ tiles that differ by a rotation by $\pi$
in the deflation of $\delta$.  In the deflation of $\kappa$ one sees 
two $\delta$ tiles that differ by reflection, and two $\delta$ tiles
that differ by rotation by $6 \pi/5$.  The reflection and these 
two rotations generate $D_{10}$.  

In the 3-dimensional dite-and-kart tiling, one has the same generators,
plus the twist of the square-faced parallelpipeds.  The twist introduces
a rotation by $\pi/2$ about a perpendicular axis, and extends the 
group to $G(10,4)$.
\vs.3
\nd {\bf \S 4.\ Algebraic and Transcendental Rotations}
\vs.1
So far we have considered groups generated by rotations by angles that
are rational multiples of $2 \pi$.  In this next example we consider
rotations by irrational multiples of $2 \pi$. Our first example comes
 from a 3-dimensional version of the pinwheel tiling [Rad] (Fig.\ 3).

Consider the right triangle $\phi$ with legs 1 and 2, and the
deflation of it given in Fig.\ 4, which decomposes $\phi$ into 25
congruent triangles each similar to $\phi$ and smaller by a linear
factor of 5. We fatten $\phi$ by width 1 to make a triangular prism we
call the ``wedge''. We now give a deflation rule for the wedge,
consisting of 5 identical layers each looking almost like Fig.\ 4 
 from one direction but
with an added complication similar to that which arose in the deflation
rules for dites. 

Note the heavily outlined rectangle in Fig.\ 4 consisting of 2 small
triangles meeting along their hypothenuses. When these triangles are
fattened to become wedges and appear in the 5 layers of the deflation
rule for a wedge, these pairs of wedges do not appear in their
original orientations but are first rotated about the axes joining the
centers of their square faces, just as we did for dites.

Given this deflation rule for wedges, ``wedge tilings'' are made by
infinite repetition of deflation-expansion-repositioning, just as for
dite and kart tilings. We now analyze the relative orientations of the
wedges in such a tiling.

Let $\nu = 2 \tan^{-1}({1 \over 2}) = \tan^{-1}({4 \over 3})$. 
In the 2 dimensional pinwheel tiling, the group of relative orientations 
is generated by rotation by $\nu$, rotation by $\pi/2$ and reflection.
In the 3 dimensional wedge tiling, the group is generated by $R_x^{\nu}$,
$R_x^{\pi/2}$ and $R_y^{\pi/2}$.  We consider the subgroup $G(\nu,4,1)$
generated by $S=R_y^{\pi/2}$ and $T=R_x^{\nu}$, and the further subgroup
$G(\nu,1,\nu)$.  Our results are extremely similar to the rational case:

\nd {\bf Lemma 6:} {\it An expression of the form
$$ W S^{b_1} T^{a_1} \cdots S^{b_n} T^{a_n} E, \eqno (4.1)  $$
with $W,E \in G(1,4,1)$, each $b_i$ odd, each $a_i$ nonzero, and $n>0$, 
cannot equal the identity matrix.}
\vs.1
\nd {\bf Theorem 4:} {\it The group $G(\nu,4,1)$ generated by $T$ and $S$
has the presentation \break
$< \a, \b : \b^4, \b^2 \a \b^2 \a >$, with the
identification $\a \to T$, $\b \to S$.}

\nd {\bf Corollary 3:} {\it The subgroup $G(\nu,1,\nu)$ of $G(\nu,4,1)$ is
isomorphic to a free group on two generators, with the generators
corresponding to $T$ and $S^{-1}TS$.}
\vs.1
\nd Remark: It is a well-known result of Stanislaw Swierczkowski that,
if $\cos(\theta)$ is rational and not equal to $0$, $\pm {1\over 2}$
or $\pm 1$, $G(\theta,1,\theta)$ is isomorphic to the free group
$<\a,\b>$, with $\a \mapsto R^\theta_x$ and $\b \mapsto R^\theta_z$
[Swi].  Since $\cos(\nu)=3/5$, Corollary 3 is a special case of
Swierczkowski's theorem.
\vs.1
\nd {\it Proof of Lemma 6}:  As in the proof of Lemma 1, we consider products
$F_a S T^a$, where $F_a$ is a numerical factor, show that all the matrix 
elements live in a certain ring $R$, and show that the (1,2), (1,3), (2,2),
and (2,3) elements (and only these elements) fail to live in a certain
maximal ideal $I$.  In this case $R = \Z$, $I$ is the
principal ideal $(5)$, $R/I = \Z_5$ and $F_a = 5^{|a|}$.

The cosine and sine of $n \nu$ are the real and imaginary parts of 
$(3 + 4i)^n/5^{n}$.  Now, if $n>0$,  the real and imaginary parts of
$(3+4i)^n$ equal 3 and 4 (mod 5), respectively.  
Thus, for any positive $a$, $5^{a} \cos(a\nu)$
and $ 5^{a} \sin(a\nu)$ are integers but not divisible by 5, while for
$a<0$, $5^{-a} \cos(a\nu)= 5^{-a} \cos(-a\nu)$
and $ 5^{-a} \sin(a\nu)= -5^{-a} \sin(a\nu)$ 
are integers but not divisible by 5.  Since
$S^{b_i}T^{a_i}$ takes the form (1.12) or (1.13), 
$F_{a_i} S^{b_i} T^{a_i}$ takes the form
$$ \pmatrix{0&\epsilon&\beta\cr 0&\gamma&\delta\cr 0&0&0\cr}\pmod{5}, 
\eqno(4.2)$$
with $\epsilon, \beta, \gamma, \delta$ nonzero elements of $\Z_5$.
But the product of two (or more) matrices of this form again takes
this form, so $FS^{b_1}T^{a_1}S^{b_2}T^{a_2}\cdots S^{b_n}T^{a_n}$,
where $F$ is the appropriate product of the $F_{a_i}$'s, again takes
this form.  Matrices in the group $G(1,4,1)$ are, up to sign,
permutation matrices, so $FWS^{b_1}T^{a_1}S^{b_2}T^{a_2}\cdots
S^{b_n}T^{a_n}E$ has 4 matrix elements that are nonzero in $\Z_5$.
But $F$ times the identity matrix is clearly zero modulo $5$, so
$WS^{b_1}T^{a_1}S^{b_2}T^{a_2}\cdots S^{b_n}T^{a_n}E$ can never equal
the identity.  \qed

\nd {\it Proof of Theorem 4}:  
The map that sends $\a \to T$ and $\b \to S$ is a well-defined
homomorphism from the abstract group to $G(\nu,4,1)$, and is clearly
onto.  We must show that it is 1-1.  Using the given relations, any
word in $\a$ and $\b$ can either be written as a power of $\b$ or as
$\b^{b_W}T^{a_1} \b T^{a_2} \cdots \b T^{a_n} \b^{b_E}$, where $n>0$
and each $a_i$ is nonzero.  By Lemma 6, the image of such an
expression in $G(\nu,4,1)$ is not the identity.  And the only powers
of $\b$ that map to the identity are powers of $\b^4=1$.  \qed

\nd {\it Proof of Corollary 3}:  Any nontrivial
word in $T$ and $S ^{-1} TS$ is of the
form (4.1), and so cannot equal the identity. \qed
\vs.1
\nd {\bf Theorem 5.}\ {\it Define the rotations $X=R_x^\omega$
and $V=S^{-1}XS=R_z^\omega$, where $x\equiv e^{i\omega}$
(equivalently $\cos(\omega)$) is transcendental. Then
the group generated by $X$ and $V$ is the free
group with those generators.}
\vs.1
\nd {\it Proof}:\ Any word in the group generated by $X$ and $V$ is of 
the form $X^{\tilde b_1}V^{\tilde d_1}X^{\tilde b_2}\cdots$ or
$V^{\tilde d_1}X^{\tilde b_1}V^{\tilde d_2}\cdots$, and can be
expressed as $X^{b_1}S^3X^{d_1}SX^{b_2}\cdots $ or
$S^3X^{d_1}SX^{b_1}S^3X^{d_2}\cdots$.  Using $S^2X^a=X^{-a}S^2$, we
can put either expression in the form
$$S^aSX^{c_1}SX^{c_2}\cdots SX^{c_n}S^b, \eqno(4.3) $$
where $a,b$ and $c_j$ are integers. All we need to show
is that $n>0$ implies that 
$S^aSX^{c_1}SX^{c_2}\cdots SX^{c_n}S^b$
is not the unit matrix. 

Each factor $SX^{c_j}$
is of the form 
$$\pmatrix{0&-\tilde s_j&\tilde c_j\cr 0&\tilde c_j&\tilde s_j\cr -1&0&0\cr},
\eqno(4.4) $$
where $\tilde c_j=\cos(c_j\omega)={x^{c_j}+x^{-c_j}\over 2},
\ \tilde s_j=\sin(c_j\omega)={x^{c_j}-x^{-c_j}\over 2i}$. 
The $(2,2)$ matrix
element of $SX^{c_1}SX^{c_2}\cdots SX^{c_n}$ is a sum of terms.
One term is the product
$\prod_j\cos(c_j \omega)=\prod_j{x^{c_j}+x^{-c_j}\over 2}$ of
the $(2,2)$ matrix elements of all the factors, and is a high-order
polynomial in $x$ and $x^{-1}$.  The remaining terms each contain at least
one power of the (3,1) element $-1$, and so are lower-order polynomials
in $x$ and $x^{-1}$.  The sum is therefore a polynomial with the same
leading term as the product $\prod_j\cos(c_j\omega)$.
Since $x$ is transcendental,
this polynomial cannot equal 0 or 1.

The factors $S^a$ and $S^b$ are, up to signs, permutations, so some
matrix element of $S^aSX^{c_1}SX^{c_2}\cdots SX^{c_n}S^b$ must be
neither 0 nor 1, and so $S^aSX^{c_1}SX^{c_2}\cdots SX^{c_n}S^b$
cannot be the unit matrix.\qed
\vs.4 \nd
{\bf Acknowledgements.}\  It is a pleasure to thank John Conway and 
Douglas Van Wieren for useful discussions; in particular, John Conway was
very helpful concerning amalgamated free products.
\vfill\eject
\nd {\bf References}
\vs.3 \nd
[CoR]\ J.\ H.\ Conway and C.\ Radin, Quaquaversal tilings and rotations,
{\it Inventiones math.}, to appear.
\vs.1 \nd
[GrS]\ B.\ Gr\"unbaum and G.C.\ Shephard, {\it Tilings and Patterns},
Freeman, New York, 1986.
\vs.1 \nd
[Rad]\ C.\ Radin, The pinwheel tilings of the plane, {\it Annals of
Math.}  {\bf 139} (1994), 661-702.
\vs.1 \nd
[Swi]\ S.\ Swierczkowski, A class of free rotation groups, {\it Indag.\ 
Math. (N.S.)} {\bf 5} (1995), 221-226.
\vfill \eject
\nopagenumbers
\hbox{}
\vs1 \nd
\hs.75 \vbox{\epsfxsize=4.5truein\epsfbox{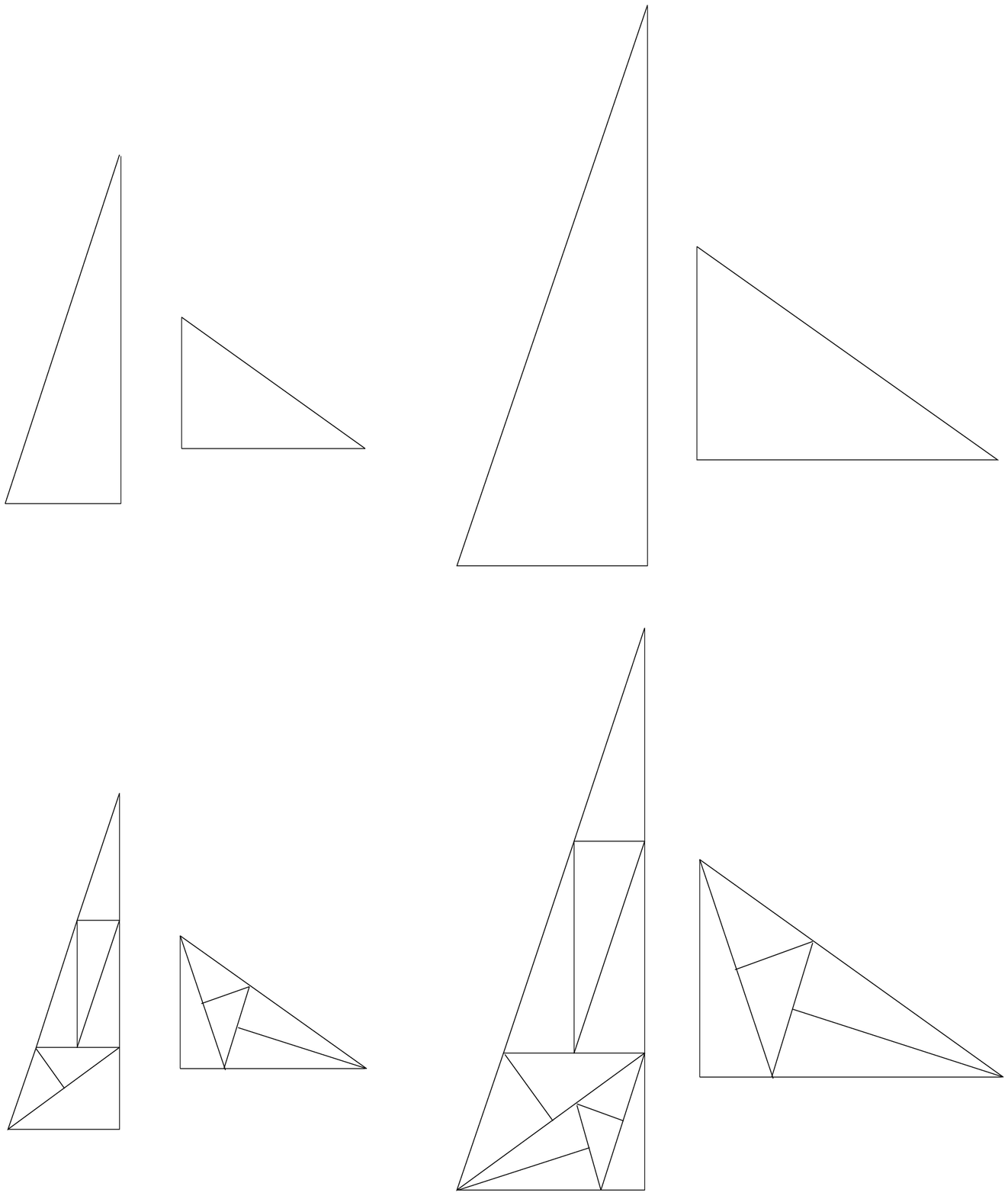}}
\vs-2.6 \nd \hs1 \hbox{$\delta$}
\hs.75 \hbox{$\kappa$}\hs1.25 \hbox{$\tilde\delta$}
\hs1.1 \hbox{$\tilde\kappa$}
\hfil
\vs2.7
\nd \hs.9 \hbox{$D(\delta)$}
\hs.4 \hbox{$D(\kappa)$}\hs1 \hbox{$D(\tilde\delta)$}
\hs.8 \hbox{$D(\tilde\kappa)$}
\hfil
\vs.5 \nd
\centerline{Figure 1. Dites and karts}
\vfill\eject
\hbox{}
\hs.75 \vbox{\epsfxsize=4.5truein\epsfbox{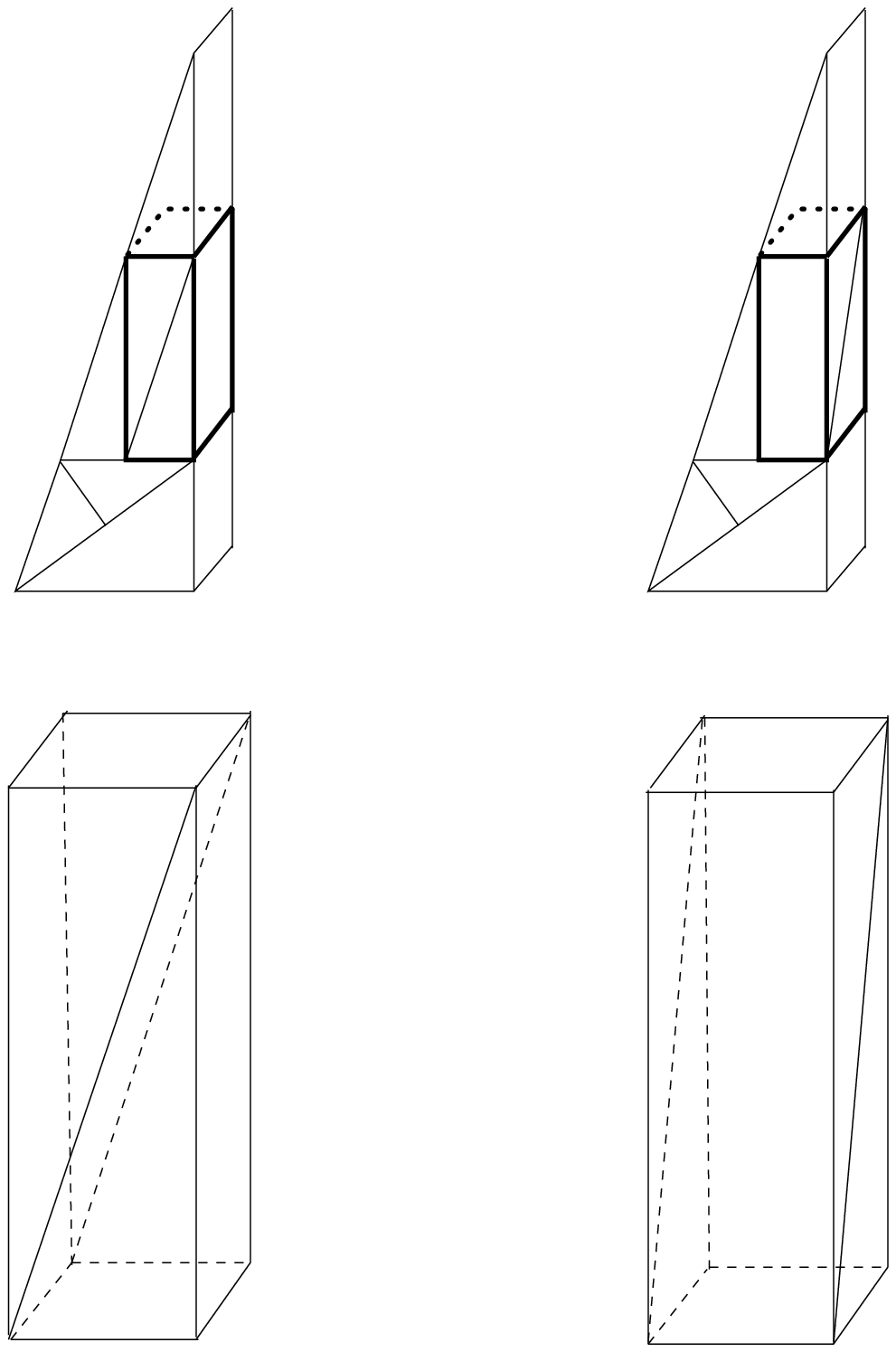}}
\vs.5
\nd \hs1.25 original
\hs2.7 rotated
\hfil
\vs.5 \nd
\centerline{Figure 2. Rotating the boxes}
\vfill\eject
\hbox{}
\vs.1 
\hbox{}\hs-1.5 \vbox{\epsfxsize=7.5truein\epsfbox{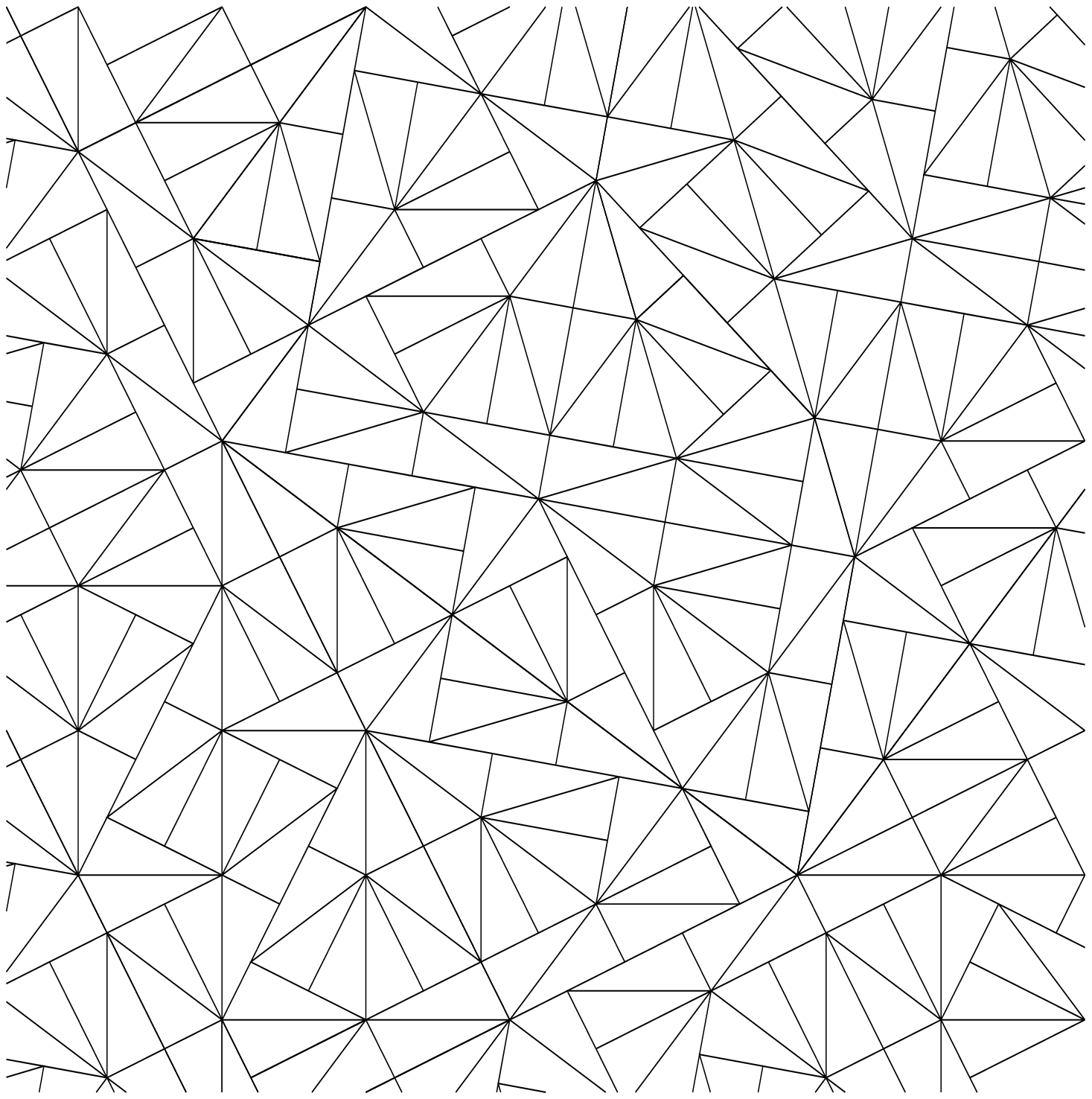}}
\vs.1 
\centerline{Figure 3. \ A pinwheel tiling}
\vfill \eject
\hbox{}
\vs2 \nd
\hs.75 \vbox{\epsfxsize=4truein\epsfbox{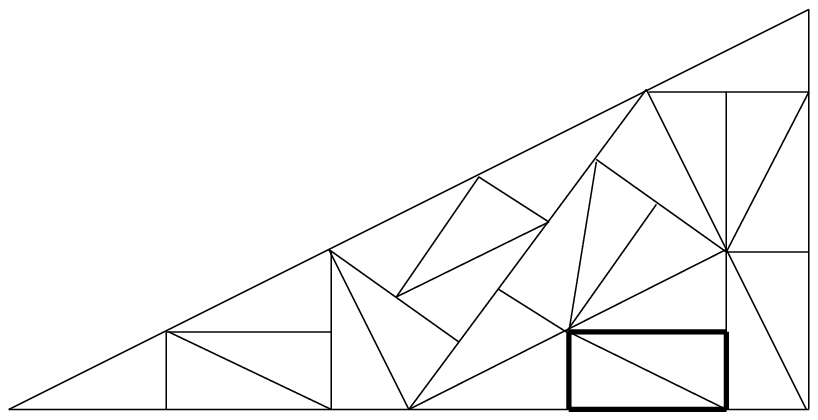}}
\vs.5 \nd
\centerline{Figure 4. Decomposition of $\phi$}
\vfill\eject

\end